\documentclass[12pt]{article}

\usepackage{amsmath,amssymb,latexsym,theorem,bbm}

\newfont{\msa}{msam10 scaled\magstep1}
\newfont{\ssmsa}{msam9}
\newfont{\smsa}{msam10}
\newfont{\sms}{msbm10}
\newfont{\sseufb}{eufb9}
\newfont{\seufb}{eufb10}
\newfont{\eufb}{eufb10 scaled\magstep1}
\newfont{\eusb}{eusb10 scaled\magstep1}
\newfont{\hcmr}{cmr17 scaled\magstep5}

\newcommand{\AstkoneKn}{\raise-7mm\hbox{\hcmr*}%
          ^{\hspace*{-5.3mm}\raise2.8mm\hbox{$\scriptstyle K_n$}}%
           _{\hspace*{-6.1mm}\raise2.5mm\hbox{$\scriptstyle k=1$}}}
\newcommand{\Astkonen}{\raise-7mm\hbox{\hcmr*}%
          ^{\hspace*{-5.3mm}\raise2.8mm\hbox{$\scriptstyle n$}}%
           _{\hspace*{-6.1mm}\raise2.5mm\hbox{$\scriptstyle k=1$}}}

\newcommand{\NN}{\mathbb{N}}

\newcommand{\RR}{\mathbb{R}}

\newcommand{\ZZ}{\mathbb{Z}}

\newcommand{\sSS}{\raise-0.5truemm\hbox{\sms S}}

\newcommand{\cM}{\mathcal{M}}

\newcommand{\cC}{\mathcal{C}}

\newcommand{\cB}{\mathcal{B}}

\newcommand{\cN}{\mathcal{N}}

\newcommand{\dd}{\mathrm{d}}

\newcommand{\sleq}{\mbox{\ssmsa\hspace*{0.1mm}\symbol{54}\hspace*{0.1mm}}}
\newcommand{\sgeq}{\mbox{\ssmsa\hspace*{0.1mm}\symbol{62}\hspace*{0.1mm}}}
\renewcommand{\leq}{\mbox{\msa\hspace*{0.9mm}\symbol{54}\hspace*{0.9mm}}}
\renewcommand{\geq}{\mbox{\msa\hspace*{0.9mm}\symbol{62}\hspace*{0.9mm}}}

\newcommand{\proofend}{\hfill\mbox{$\Box$}}

\numberwithin{equation}{section}

\theoremstyle{change} \theorembodyfont{\em}
\newtheorem{Lem}{Lemma.}[section]
\newtheorem{Thm}[Lem]{Theorem.}

\theorembodyfont{\rm}

\newtheorem{Rem}[Lem]{Remark.}

\begin{document}

\begin{center}
  {\bfseries\Large Portmanteau theorem \\
         for unbounded measures}\\ [5mm]
  By\\[3mm]
  {\bfseries\large M\'aty\'as Barczy} {\large and}
  {\bfseries\large Gyula Pap}\\[3mm]
  University of Debrecen, Hungary
\end{center}

\renewcommand{\thefootnote}{}
\footnote{The authors have been supported by the Hungarian
Scientific Research Fund under Grant No.\ OTKA--T048544/2005. The
first author has been also supported by the Hungarian Scientific
Research Fund under Grant No.\ OTKA--F046061/2004.}

\vspace*{-2mm}

{\small\textbf{Abstract.} We prove an analogue of the portmanteau
theorem on weak convergence of probability measures allowing
measures which are unbounded on an underlying metric space but
finite on the complement of any Borel neighbourhood of a fixed
element.


\section{Introduction}

Weak convergence of probability measures on a metric space has a
very important role in probability theory. The well known
\emph{portmanteau theorem} due to A. D. Alexandroff (see for
 example Theorem 11.1.1 in Dudley \cite{DUD}) provides useful conditions
 equivalent to weak convergence of probability measures; any of them could
 serve as the definition of weak convergence.
Proposition 1.2.13 in the book of Meerschaert and Scheffler
\cite{MS} gives an analogue of the portmanteau theorem for bounded
measures on \ $\RR^d$. \ Moreover, Proposition 1.2.19 in \cite{MS}
gives an analogue for special unbounded measures on \ $\RR^d$, \
more precisely, for extended real valued measures which are finite
on the complement of any Borel neighbourhood of \ $0\in\RR^d$.

By giving counterexamples we show that the equivalences of \ (c) \
and \ (d) \ in Propositions 1.2.13 and 1.2.19 in \cite{MS} are not
valid (see our Remarks \ref{REMARK1} and \ref{REMARK2}). We
reformulate Proposition 1.2.19 in \cite{MS} in a more detailed form
adding new equivalent assertions to it (see Theorem
\ref{TETEL:PORTMANTEAU}). Moreover, we note that Theorem
\ref{TETEL:PORTMANTEAU} generalizes the equivalence of \ (a) \ and \
(b) \ in Theorem 11.3.3 of \cite{DUD} in two aspects. On the one
hand, the equivalence is extended allowing not necessarily finite
 measures which are finite on the complement of any Borel neighbourhood of a
 fixed element of an underlying metric space.
On the other hand, we do not assume the separability of the
underlying metric space to prove the equivalence. But we mention
that this latter possibility is hiddenly contained in Problem 3,
 p. 312 in \cite{DUD}. For completeness we give a detailed proof of Theorem
\ref{TETEL:PORTMANTEAU}. Our proof goes along the lines of the proof
of the original portmanteau theorem and differs from the proof of
Proposition 1.2.19 in \cite{MS}.

To shed some light on the sense of a portmanteau theorem for
unbounded measures, let us consider the question of weak convergence
of infinitely divisible probability measures \ $\mu_n$, \ $n\in\NN$
\ towards an infinitely divisible probability measure \ $\mu_0$ \ in
case of the real line \ $\RR$. \ Theorem VII.2.9 in Jacod and
Shiryayev \cite{JSH} gives equivalent conditions
 for weak convergence \ $\mu_n\stackrel{w}{\rightarrow}\mu_0$.
\ Among these conditions we have
 \begin{align}\label{LEVY_KONVERG_FELTETEL}
  \int_{\RR}f\,\dd\eta_n\to\int_{\RR}f\,\dd\eta_0
  \qquad\text{for all \ $f\in\cC_2(\RR)$,}
 \end{align}
 where \ $\eta_n$, \ $n\in\ZZ_+$ \ are nonnegative, extended real valued
 measures on \ $\RR$ \ with \ $\eta_n(\{0\})=0$ \ and
 \ $\int_{\RR}(x^2\wedge 1)\,\dd\eta_n(x)<\infty$
 \ (i.e., L\'evy measures on \ $\RR$) \ corresponding to \ $\mu_n$, \ and
 \ $\cC_2(\RR)$ \ is the set of all real valued bounded continuous functions
 \ $f$ \ on \ $\RR$ \ vanishing on some Borel neighbourhood of \ $0$ \ and
 having a limit at infinity.
Theorem \ref{TETEL:PORTMANTEAU} is about equivalent reformulations
of \eqref{LEVY_KONVERG_FELTETEL} when it holds for all real valued
bounded continuous functions on \ $\RR$ \ vanishing on some Borel
neighbourhood of \ $0$.

\section{An analogue of the portmanteau theorem}

Let \ $\NN$ \ and \ $\ZZ_+$ \ be the set of positive and nonnegative
integers, respectively. Let \ $(X,d)$ \ be a metric space and \
$x_0$ \ be a fixed element of \ $X$. \ Let \ $\cB(X)$ \ denote the
$\sigma$-algebra of Borel subsets of \ $X$. \ A Borel neighbourhood
\ $U$ \ of \ $x_0$ \ is an element of \ $\cB(X)$ \ for
 which there exists an open subset \ $\widetilde U$ \ of \ $X$ \ such that
 \ $x_0\in\widetilde{U}\subset U$.
\ Let \ $\cN_{x_0}$ \ denote the set of all Borel neighbourhoods of
\ $x_0$, \ and the set of bounded measures on \ $X$ \ is denoted by
\ $\cM^b(X)$. \ The expression ''a measure \ $\mu$ \ on \ $X$''
means a measure \ $\mu$ \ on the $\sigma$-algebra \ $\cB(X)$.

 Let \ $\cC(X)$, \ $\cC_{x_0}(X)$ \ and \ ${\rm BL}_{x_0}(X)$ \ denote the
 spaces of all real valued bounded continuous functions on \ $X$, \ the set of
 all elements of \ $\cC(X)$ \ vanishing on some Borel neighbourhood of \ $x_0$,
 \ and the set of all real valued bounded
Lipschitz functions vanishing on some Borel neighbourhood of \
$x_0$, \ respectively.

For a measure \ $\eta$ \ on \ $X$ \ and for a Borel subset \
$B\in\cB(X)$,
 \ let \ $\eta|_B$ \ denote the restriction of \ $\eta$ \ onto \ $B$, \ i.e.,
 \ $\eta|_B(A):=\eta(B\cap A)$ \ for all \ $A\in\cB(X)$.

Let \ $\mu_n$, \ $n\in\ZZ_+$ \ be bounded measures on \ $X$.
 \ We write \ $\mu_n\stackrel{w}{\rightarrow}\mu$ \ if \ $\mu_n(A)\to\mu(A)$
 \ for all \ $A\in\cB(X)$ \ with \ $\mu(\partial A)=0$.
\ This is called {\em weak convergence of bounded measures} on \
$X$.

Now we formulate a portmanteau theorem for unbounded measures.

\begin{Thm}\label{TETEL:PORTMANTEAU}
 Let \ $(X,d)$ \ be a metric space and \ $x_0$ \ be a fixed element of \ $X$.
\ Let \ $\eta_n$, \ $n\in\mathbb Z_+$, \ be measures on \ $X$ \ such
that
 \ $\eta_n(X\setminus U)<\infty$ \ for all \ $U\in\cN_{x_0}$ \ and for all
 \ $n\in\mathbb Z_+$.
\ Then the following assertions are equivalent:
 \renewcommand{\labelenumi}{{\rm(\roman{enumi})}}
 \begin{enumerate}
  \item $\int_{X\setminus U}f\,\dd\eta_n\to\int_{X\setminus U}f\,\dd\eta_0$
         \ for all \ $f\in\cC(X)$, \ $U\in\cN_{x_0}$ \ with
         \ $\eta_0(\partial U)=0$,
  \item \ $\eta_n|_{X\setminus U}
           \stackrel{w}{\rightarrow}\eta_0|_{X\setminus U}$
        \ for all \ $U\in\cN_{x_0}$ \ with \ $\eta_0(\partial U)=0$,
  \item \ $\eta_n(X\setminus U)\to\eta_0(X\setminus U)$
         \ for all \ $U\in\cN_{x_0}$ \ with \ $\eta_0(\partial U)=0$,
  \item $\int_X f\,\dd\eta_n\to\int_X f\,\dd\eta_0$ \ for all
         \ $f\in\cC_{x_0}(X)$,
  \item $\int_Xf\,\dd\eta_n\to\int_Xf\,\dd\eta_0$ \ for all
         \ $f\in {\rm BL}_{x_0}(X)$,
  \item the following inequalities hold:\\
   \noindent \hspace*{-0.5cm}
    {\rm(}$a${\rm)}
    \ $\limsup\limits_{n\to\infty}\eta_n(X\setminus U)\leq\eta_0(X\setminus U)$
    \ for all open neighbourhoods \ $U$ \ of \ $x_0$,\\
   \noindent \hspace*{-0.5cm}
    {\rm(}$b${\rm)} \
    \ $\liminf\limits_{n\to\infty}\eta_n(X\setminus V)\geq\eta_0(X\setminus V)$
    \ for all closed neighbourhoods \ $V$ \ of \ $x_0$.
 \end{enumerate}
\end{Thm}

\noindent \textbf{Proof.} {\bf (i)$\Rightarrow$(ii):} \ Let \ $U$ \
be an element of \ $\cN_{x_0}$ \ with \ $\eta_0(\partial U)=0$. \
Note \ $\eta_n|_{X\setminus U}\in\cM^b(X),$ $n\in\ZZ_+$. \ By the
equivalence of \ (a) \ and \ (b) \ in Proposition 1.2.13 in
\cite{MS}, to prove
 \ $\eta_n|_{X\setminus U}\stackrel{w}{\rightarrow}\eta_0|_{X\setminus U}$
 \ it is enough to check
 \ $\int_X f\,\dd\eta_n|_{X\setminus U}\to\int_X f\,\dd\eta_0|_{X\setminus U}$
 \ for all \ $f\in\cC(X)$.
\ For this it suffices to show that for all real valued bounded
measurable functions \ $h$ \ on \ $X$, \ for all \ $A\in\cB(X)$ \
and for all \ $n\in\ZZ_+$ \ we have
 \begin{align}\label{SEGED2}
  \int_Xh\,\dd\eta_n|_{A}
  =\int_{A}h\,\dd\eta_n.
 \end{align}
By Beppo-Levi's theorem, a standard measure-theoretic argument
implies \eqref{SEGED2}.

\noindent {\bf (ii)$\Rightarrow$(iii):} \ Let \ $U$ \ be an element
of \ $\cN_{x_0}$ \ with \ $\eta_0(\partial U)=0$. \ By \ (ii), we
have
 \ $\eta_n|_{X\setminus U}\stackrel{w}{\rightarrow}\eta_0|_{X\setminus U}$.
\ Since
 \ $\eta_0|_{X\setminus U}(\partial X)=\eta_0|_{X\setminus U}(\emptyset)=0$,
 \ we get
 \ $\eta_n(X\setminus U)=\eta_n|_{X\setminus U}(X)
    \to\eta_0|_{X\setminus U}(X)=\eta_0(X\setminus U)$,
 \ as desired.

\noindent {\bf (iii)$\Rightarrow$(ii):} \ Let \ $U$ \ be an element
of \ $\cN_{x_0}$ \ with \ $\eta_0(\partial U)=0$
 \ and let \ $B\in\cB(X)$ \ be such that
 \ $\eta_0|_{X\setminus U}(\partial B)=0$.
\ We have to show
 \ $\eta_n|_{X\setminus U}(B)\to\eta_0|_{X\setminus U}(B)$.

Since \ $B\cap(X\setminus U)=X\setminus[X\setminus(B\cap(X\setminus
U))]$ \ and \ $\eta_n|_{X\setminus U}(B)=\eta_n(B\cap(X\setminus
U))$, \ $n\in\ZZ_+$, \ by \ (iii), \ it is enough to check
 \ $\eta_0\big(\partial\big(X\setminus(B\cap(X\setminus U))\big)\big)=0$.
\ First we show
 \begin{align}\label{SEGED1}
  \text{$\partial\big(B\cap(X \! \setminus \! U)\big)
         \subset\big(\partial B\cap(X \! \setminus \! U)\big)\cup\partial U$
        \ for all subsets \ $B$, $U$ \ of \ $X$.}
 \end{align}
Let \ $x$ \ be an element of \ $\partial\big(B\cap(X\setminus
U)\big)$ \ and
 \ $(y_n)_{n\sgeq 1}$, \ $(z_n)_{n\sgeq 1}$ \ be two sequences such that
 \ $\lim_{n\to\infty}y_n=\lim_{n\to\infty}z_n=x$ \ and
 \ $y_n\in B\cap(X\setminus U)$, \ $z_n\in X\setminus(B\cap(X\setminus U))$,
 \ $n\in\mathbb N$.
\ Then for all \ $n\in\mathbb N$ \ we have one or two of the
following possibilities:
 \begin{itemize}
   \item $y_n\in B,$ \ $y_n\in X\setminus U$ \ and \ $z_n\in X\setminus B,$
   \item $y_n\in B,$ \ $y_n\in X\setminus U$ \ and \ $z_n\in U.$
 \end{itemize}
Then we get
 \ $x\in\big(\partial B\cap((X\setminus U)\cup\partial U)\big)
         \cup\big(\partial U\cap (B\cup\partial B)\big)
         \cup\big(\partial B\cap \partial U\big)
   $.
\ Since
 \ $\partial B\cap((X\setminus U)\cup\partial U)
        \subset(\partial B\cap(X\setminus U))\cup\partial U$,
 \ we have \ $x\in\big(\partial B\cap(X\setminus U)\big)\cup\partial U$,
 \ as desired.

Using \eqref{SEGED1} we get
 \ $\eta_0\big(\partial\big(X\setminus(B\cap(X\setminus U))\big)\big)
     \leq \eta_0\big(\partial B\cap(X\setminus U)\big)
         +\eta_0(\partial U)=0$.
\ Indeed, by the assumptions
 \ $\eta_0\big(\partial B\cap(X\setminus U)\big)=0$ \ and
 \ $\eta_0(\partial U)=0$.
\ Hence
  \ $\eta_0\big(\partial\big(X\setminus(B\cap(X\setminus U))\big)\big)=0$.

\noindent {\bf (ii)$\Rightarrow$(i):} \ Using again the equivalence
of \ (a) \ and \ (b) \ in Proposition 1.2.13 in
 \cite{MS} and \eqref{SEGED2} \ we obtain (i).

\noindent {\bf (iii)$\Rightarrow$(iv):} \ Let \ $f$ \ be an element
of \ $\cC_{x_0}(X)$. \ Then there exists \ $A\in\cN_{x_0}$ \ such
that \ $f(x)=0$ \ for all
 \ $x\in A$ \ and \ $\eta_0(\partial A)=0$.
\ Indeed, the function \ $t\mapsto\eta_0\big(\{x\in X:d(x,x_0)\geq
t\}\big)$
 \ from \ $(0,+\infty)$ \ into \ $\RR$ \ is monotone decreasing, hence the set
 \ $\big\{t\in(0,+\infty):\eta_0(\{x\in X:d(x,x_0)=t\})>0\big\}$ \ of its
 discontinuities is at most countable.
Consequently, for all \ $\widetilde U\in\cN_{x_0}$ \ there exists
some \ $t>0$
 \ such that \ $U:=\{x\in X:d(x,x_0)<t\}\in\cN_{x_0}$, \ $U\subset\widetilde U$
 \ and \ $\eta_0(\partial U)=0$.
\ (At this step we use that an element \ $\widetilde U$ \ of \
$\cN_{x_0}$ \ contains an open subset of \ $X$ \ containing \
$x_0$.) \ This implies the existence of \ $A$. \ We show that the
set
 \ $
     D:=\big\{t\in\mathbb R:\eta_0\big(\{x\in X:f(x)=t\}\big)>0\big\}
   $
 \ is at most countable.
The function \ $F:\mathbb R\to[0,\eta_0(X\setminus A)]$, \ defined
by
 $$
  F(t):=\eta_0\big(\{x\in X\setminus A:f(x)<t\}\big),
            \quad t\in\mathbb R,
 $$
 is monotone increasing and left continuous.
(Note that \ $\eta_0(X\setminus A)<\infty,$ \ by the assumption on \
$\eta_0.$) \ Hence it has at most countably many discontinuity
points, and
 \ $t_0\in\mathbb R$ \ is a discontinuity point of \ $F$ \ if and only if
 \ $F(t_0+0)>F(t_0),$ \ i.e.,
  \ $
     \eta_0\big(\{x\in X\setminus A: f(x)=t_0\}\big)>0
    $.
\ If \ $t_0\ne 0,$ \ then
 \ $
    \{x\in X: f(x)=t_0\}=\{x\in X\setminus A: f(x)=t_0\}
   $,
\ thus \ $t_0\ne 0$ \ is a discontinuity point of \ $F$ \ if and
only if
 \ $\eta_0(\{x\in X: f(x)=t_0\})>0$.
\ Hence if \ $t\in D$ \ then \ $t=0$ \ or \ $t$ \ is a discontinuity
point of
 \ $F$, \ consequently \ $D$ \ is at most countable.
Since \ $f$ \ is bounded and \ $D$ \ is at most countable, there
exists a real
 number \ $M>0$ \ such that \ $-M,M\notin D$ \ and \ $\vert f(x)\vert<M$ \ for
 \ $x\in X$.
\ Let \ $\varepsilon>0$. \ Choose real numbers \ $t_i,$
$i=0,\ldots,k$ \ such that
 \ $-M=t_0<t_1<\cdots<t_k=M$, \ $t_i\notin D,$ $i=0,\ldots,k$ \ and
 \ $\max_{0\sleq i\sleq k-1}(t_{i+1}-t_i)<\varepsilon$.
\ The countability of \ $D$ \ implies the existence of \ $t_i$,
 \ $i=0,\ldots,k$.
\ Let
 $$
  B_i:=f^{-1}\big([t_i,t_{i+1})\big)\cap(X\setminus A)
      =\Big\{x\in X\setminus A:t_i\leq f(x)<t_{i+1}\Big\}
 $$
 for all \ $i=0,\ldots,k-1$.
\ Then \ $B_i$, \ $i=0,\ldots,k-1$, \ are pairwise disjoint Borel
sets and
  \ $X\setminus A=\bigcup_{i=0}^{k-1}B_i$.
\ Since \ $f$ \ is continuous, the boundary \ $\partial(f^{-1}(H))$
\ of the
 set \ $f^{-1}(H)$ \ is a subset of the set \ $f^{-1}(\partial H)$ \ for all
 subsets \ $H$ \ of \ $\RR$.
\ Using \eqref{SEGED1} this implies
 \ $
    \partial(X\setminus B_i)=\partial B_i
               \subset f^{-1}(\{t_i\})\cup f^{-1}(\{t_{i+1}\})
                     \cup\partial A
   $
 \ for all \ $i=0,\ldots,k-1$.
\ Since \ $t_i\notin D,$ $i=0,\ldots,k,$ \ $\eta_0(\partial A)=0$
 \ and
 $$
  \eta_0(\partial(X \! \setminus \! B_i))
     \leq\eta_0\big(\{x \! \in \! X:f(x) \! = \! t_i\}\big)
         +\eta_0\big(\{x \! \in \! X:f(x) \! = \! t_{i+1}\}\big)
         +\eta_0(\partial A),
 $$
 we get \ $\eta_0(\partial(X\setminus B_i))=0$, \ $i=0,\ldots,k-1$.
\ Since \ $A\subset X\setminus B_i$, \ we have \ $X\setminus
B_i\in\cN_{x_0}$
 \ for all \ $i=0,\ldots,k-1$.
\ Hence condition \ (iii) \ implies that \
$\eta_n(B_i)\to\eta_0(B_i)$ \ as
 \ $n\to\infty,$ $i=0,\ldots,k-1$.
\ By the triangle inequality
 \begin{align*}
  \Big\vert\int_Xf\,\dd\eta_n-\int_Xf\,\dd\eta_0\Big\vert
    &=\Big\vert\int_{X\setminus A}f\,\dd\eta_n
               -\int_{X\setminus A}f\,\dd\eta_0\Big\vert\\
    &\,\leq 2\max_{0\sleq i\sleq k-1}(t_{i+1}-t_i)
        +\Big\vert\sum_{i=0}^{k-1}
                   t_i\big(\eta_n(B_i)-\eta_0(B_i)\big)\Big\vert.
 \end{align*}
Hence
 \ $
    \limsup_{n\to\infty}
     \left\vert\int_Xf\,\dd\eta_n-\int_Xf\,\dd\eta_0\right\vert
    \leq 2\max_{0\sleq i\sleq k-1}(t_{i+1}-t_i)
    <2\varepsilon
  $.
\ Since \ $\varepsilon>0$ \ is arbitrary, \ (iv) \ holds.

\noindent {\bf (iv)$\Rightarrow$(v):} \ It is trivial, since \ ${\rm
BL}_{x_0}(X)\subset\cC_{x_0}(X)$.

\noindent {\bf (v)$\Rightarrow$(vi):} \ First let \ $U$ \ be an open
neighbourhood of \ $x_0$. \ Let \ $\varepsilon>0$. \ We show the
existence of a closed neighbourhood \ $U_\varepsilon$ \ of
 \ $x_0$ \ such that \ $U_\varepsilon\subset U$ \ and
 \ $\eta_0(U\setminus U_\varepsilon)<\varepsilon$, \ and of a function
 \ $f\in {\rm BL}_{x_0}(X)$ \ such that \ $f(x)=0$ \ for
 \ $x\in U_\varepsilon$, \ $f(x)=1$ \ for \ $x\in X\setminus U$ \ and
 \ $0\leq f(x)\leq 1$ \ for \ $x\in X$.

For all \ $B\in\cB(X)$ \ and for all \ $\lambda>0$ \ we use notation
 \ $
    B^{\lambda}:=\big\{x\in X:d(x,B)<\lambda\big\}
   $,
\ where \ $d(x,B):=\inf\{d(x,z):z\in B\}$. \ Since \ $U$ \ is open,
we get \ $U=\bigcup_{n=1}^{\infty}F_n$, \ where
 \ $F_n:=X\setminus (X\setminus U)^{1/n}$, \ $n\in\mathbb N$.
\ Then \ $F_n\subset F_{n+1}$, \ $n\in\mathbb N$, \ $F_n$ \ is a
closed subset
 of \ $X$ \ for all \ $n\in\mathbb N$ \ and
 \ $\bigcap_{n=1}^{\infty}(X\setminus F_n)=X\setminus U$.
\ We also have \ $\eta_0(X\setminus F_N)<\infty$ \ for some
sufficiently large
 \ $N\in\NN$ \ and \ $X\setminus F_n\supset X\setminus F_{n+1}$ \ for all
 \ $n\in\mathbb N$, \ and hence the continuity of the measure \ $\eta_0$
 \ implies that
 \ $\lim_{n\to\infty}\eta_0(X\setminus F_n)=\eta_0(X\setminus U)$.
\ Since \ $\eta_0(X\setminus U)<\infty,$ \ there exists some
 \ $n_0\in\mathbb N$ \ such that
 \ $\eta_0(X\setminus F_{n_0})-\eta_0(X\setminus U)<\varepsilon$.
\ Set \ $U_\varepsilon:= F_{n_0}$. \ Since
 \ $
    \eta_0(X\setminus F_{n_0})-\eta_0(X\setminus U)
         =\eta_0\big((X\setminus F_{n_0})\setminus (X\setminus U)\big)
         =\eta_0(U\setminus F_{n_0})
   $,
 \ the set \ $U_\varepsilon$ \ is a closed neighborhood of \ $x_0$,
 \ $U_\varepsilon\subset U$ \ and
 \ $\eta_0(U\setminus U_\varepsilon)<\varepsilon$.

We show that the function \ $f:X\to\RR,$ \ defined by
 \ $ f(x):=\min(1,n_0d(x,U_\varepsilon))$, \ $x\in X$, \ is an element of
 \ ${\rm BL}_{x_0}(X),$ \ $f(x)=0$ \ for \ $x\in U_\varepsilon$,
 \ $f(x)=1$ \ for \ $x\in X\setminus U$ \ and \ $0\leq f(x)\leq 1$ \ for
 \ $x\in X$.

\noindent Note that if \ $x\in U_\varepsilon$ \ then \
$d(x,U_\varepsilon)=0$,
 \ hence \ $f(x)=0$.
\ And if \ $x\in X\setminus U$ \ then
 \ $d(x,U_\varepsilon)\geq d(X\setminus U,U_\varepsilon)\geq1/n_0$,
 \ hence \ $f(x)=1$.
\ The fact that \ $0\leq f(x)\leq 1,$ \ $x\in X$ \ is obvious. To
prove that \ $f$ \ is Lipschitz, we check that
 $$
  \vert f(x)-f(y)\vert\leq n_0d(x,y)
    \qquad\text{for all \ $x,y\in X$.}
 $$
If \ $x,y\in X$ \ with \ $d(x,y)\geq 1/n_0$ \ then
  \ $\vert f(x)-f(y)\vert\leq 1\leq n_0d(x,y)$.
\ If \ $x,y\in X$ \ with \ $d(x,y)<1/n_0$ \ then we have to consider
the
 following four cases apart from changing the role of \ $x$ \ and  \ $y$:
 \ $x\in X\setminus U$, \ $y\in U\setminus U_\varepsilon$;
 \ $x\in U_\varepsilon$, \ $y\in U\setminus U_\varepsilon$;
 \ $x,y\in U\setminus U_\varepsilon$ \ and the case \ $x,y\in U_\varepsilon$
 \ or \ $x,y\in X\setminus U$.

\noindent Let us consider the case when
 \ $x,y\in U\setminus U_\varepsilon$ \ and \ $f(x)=1$,
 \ $f(y)=n_0d(y,U_\varepsilon)$.
\ Then
 \ $d(x,U_\varepsilon)\geq 1/n_0,$ \ $d(y,U_\varepsilon)\leq 1/n_0$
 \ and we get \ $\vert f(x)-f(y)\vert=1-n_0d(y,U_\varepsilon)\leq n_0d(x,y)$.
\ Indeed, \ $1/n_0\leq d(x,U_\varepsilon)\leq
d(x,y)+d(y,U_\varepsilon)$. \ The case \ $x,y\in U\setminus
U_\varepsilon$ \ and \ $f(y)=1$
 \ $f(x)=n_0d(x,U_\varepsilon)$ \ can be handled similarly.
If \ $x,y\in U\setminus U_\varepsilon$ \ and
 \ $f(x)=n_0d(x,U_\varepsilon)$, \ $f(y)=n_0d(y,U_\varepsilon)$ \ then
 $$
   \vert f(x)-f(y)\vert=n_0\vert d(x,U_\varepsilon)-d(y,U_\varepsilon)\vert
        \leq n_0d(x,y).
 $$
Indeed, since \ $U_\varepsilon$ \ is closed, we have
 \ $\vert d(x,U_\varepsilon)-d(y,U_\varepsilon)\vert\leq d(x,y)$.
\ If \ $x,y\in U\setminus U_\varepsilon$ \ and \ $f(x)=f(y)=1$ \
then
 \ $\vert f(x)-f(y)\vert=0\leq n_0d(x,y)$.

\noindent The other cases can be handled similarly. Hence \ $f\in
{\rm BL}_{x_0}(X)$ \ and we get
 \begin{align*}
  &\int_X\!f\,\dd\eta_0=
                 \int_{X\setminus U_\varepsilon}\!\!f\,\dd\eta_0
                 \leq\eta_0(X \! \setminus \! U_\varepsilon)
                 =\eta_0(X \! \setminus \! U)
                  +\eta_0(U \! \setminus \! U_\varepsilon)
                 <\eta_0(X \! \setminus \! U)\!+\!\varepsilon,
 \end{align*}
 and
 \ $\int_Xf\,\dd\eta_n\geq\int_{X\setminus U}f\,\dd\eta_n
     =\eta_n(X\setminus U)$.
\ Hence by condition \ (v) \ we have
 \begin{align*}
   \limsup_{n\to\infty}\eta_n(X\setminus U)
          \leq\limsup_{n\to\infty}\int_Xf\,\dd\eta_n
          =\int_Xf\,\dd\eta_0
          <\eta_0(X\setminus U)+\varepsilon.
 \end{align*}
Since \ $\varepsilon>0$ \ is arbitrary, we get \ $(a)$.

\noindent Now let \ $V$ \ be a closed neighbourhood of \ $x_0$. \
Let \ $\varepsilon>0$. \ As in case of an open neighbourhood of \
$x_0$, \ one can show that there
 exist an open neighbourhood \ $V_\varepsilon$ \ of \ $x_0$ \ such that
 \ $V\subset V_\varepsilon$ \ and
 \ $\eta_0(V_\varepsilon\setminus V)<\varepsilon$ \ and
 a function \ $f\in {\rm BL}_{x_0}(X)$ \ such that \ $f(x)=0$ \ for \ $x\in V$,
 \ $f(x)=1$ \ for \ $x\in X\setminus V_\varepsilon$ \ and \ $0\leq f(x)\leq 1$
 \ for \ $x\in X$.
\ Then we get
 \begin{align*}
 \int_Xf\,\dd\eta_0&=\int_{X\setminus V}f\,\dd\eta_0
                   =\eta_0(X\setminus V_\varepsilon)
                            +\int_{V_\varepsilon\setminus
                            V}f\,\dd\eta_0\\
                   &\geq\eta_0(X\setminus V)-\eta_0(V_\varepsilon\setminus V)
                    >\eta_0(X\setminus V)-\varepsilon,
 \end{align*}
 and
 \ $\int_Xf\,\dd\eta_n=\int_{X\setminus V}f\,\dd\eta_n
                    \leq\eta_n(X\setminus V)$.
\ Hence by condition \ (v) \ we have
 \begin{align*}
   \liminf_{n\to\infty}\eta_n(X\setminus V)
          \geq\liminf_{n\to\infty}\int_Xf\,\dd\eta_n
          =\int_Xf\,\dd\eta_0
          >\eta_0(X\setminus V)-\varepsilon.
 \end{align*}
Since \ $\varepsilon>0$ \ is arbitrary, we obtain \ $(b)$.

{\bf (vi)$\Rightarrow$ (iii):} \ The proof can be carried out
similarly to the proof of the corresponding part
 of Theorem 11.1.1 in Dudley \cite{DUD}.
\proofend

\begin{Rem}
 Assertion \ (v) \ in Theorem \ref{TETEL:PORTMANTEAU} can be replaced by
  $$
   \int_X f\,\dd\eta_n\to\int_X f\,\dd\eta_0
    \qquad\text{for all \ $f\in\cC_{x_0}^u(X)$,}
  $$
 where \ $\cC_{x_0}^u(X)$ \ denotes the set of all uniformly continuous
 functions in $\cC_{x_0}(X)$.
\end{Rem}

\begin{Rem}\label{REMARK1}
By giving a counterexample we show that \ (a) \ and \ (b) \ in
condition \ (vi)
 \ of Theorem \ref{TETEL:PORTMANTEAU} are not equivalent.
For all \ $n\in\NN$ \ let \ $\eta_n$ \ be the Dirac measure \
$\delta_{2}$ \ on
 \ $\RR$ \ concentrated on \ $2$ \ and let \ $\eta_0$ \ be the Dirac measure
 \ $\delta_0$ \ on \ $\RR$ \ concentrated on \ 0.
Then \ $\eta_0(\RR\setminus V)=0$ \ for all closed neighbourhoods \
$V$ \ of \ 0, \ hence \ (b) \ in condition \ (vi) \ of Theorem
\ref{TETEL:PORTMANTEAU} holds. But \ (a) \ in condition \ (vi) \ of
Theorem \ref{TETEL:PORTMANTEAU} is not satisfied. Indeed, \
$U:=(-1,1)$ \ is an open neighbourhood of \ 0,
 \ $\eta_0(\RR\setminus U)=0$, \ but
 $$
  \eta_n(\RR\setminus U)
      =\eta_n\big((-\infty,-1],[1,+\infty)\big)
      =1,\quad n\in\NN,
 $$
 hence \ $\limsup_{n\to\infty}\eta_n(\RR\setminus U)=1$.
\ This counterexample also implies that the equivalence of \ $(c)$ \
and
 \ $(d)$ \ in Proposition 1.2.19 in \cite{MS} is not valid.
\end{Rem}

\begin{Rem}\label{REMARK2}
By giving a counterexample we show that the equivalence of \ (c) \
and \ (d) \ in Proposition 1.2.13 in \cite{MS} is not valid. For all
\ $n\in\NN$ \ let \ $\mu_n$ \ be the measure \ $2\delta_{1/n}$ \ on
 \ $\RR$ \ and \ $\mu$ \ be the Dirac measure \ $\delta_0$ \ on \ $\RR$.
\ We have \ $\mu(A)\leq\liminf_{n\to\infty}\mu_n(A)$ \ for all open
subsets \ $A$ \ of \ $\RR$ \ but there exists some closed subset \
$F$ \ of \ $\RR$ \ such that \
$\limsup_{n\to\infty}\mu_n(F)>\mu(F)$. \ If \ $A$ \ is an open
subset of \ $\RR$ \ such that \ $0\in A$ \ then
 \ $\mu(A)=1$ \ and \ $\mu_n(A)=2$ \ for all sufficiently large \ $n$, \ which
 implies \ $\mu(A)\leq\liminf_{n\to\infty}\mu_n(A)$.
\ If \ $A$ \ is an open subset of \ $\RR$ \ such that \ $0\notin A$
\ then \ $\mu(A)=0$, \ hence \
$\mu(A)\leq\liminf_{n\to\infty}\mu_n(A)$ \ is valid. Let \ $F$ \ be
the closed interval \ $[-1,1]$. \ Then \ $\mu(F)=1$ \ and \
$\mu_n(F)=2$, \ $n\in\NN$, \ which yields
 \ $\limsup_{n\to\infty}\mu_n(F)=2$.
\ Hence \ $\limsup_{n\to\infty}\mu_n(F)>\mu(F)$.
\end{Rem}

\parbox{42mm}{M\'aty\'as Barczy\\
             Faculty of Informatics\\
             University of Debrecen\\
             Pf.12\\
             H--4010 Debrecen\\
             Hungary\\ \\
             barczy@inf.unideb.hu}
\hfill
\parbox{42mm}{Gyula Pap\\
             Faculty of Informatics\\
             University of Debrecen\\
             Pf.12\\
             H--4010 Debrecen\\
             Hungary\\ \\
             papgy@inf.unideb.hu}

\end{document}